\documentclass[a4, 12pt,draft]{amsart}
\usepackage[mathscr]{eucal}
\usepackage{amssymb}
\usepackage{latexsym}
\usepackage{amsthm}

\theoremstyle{plain}
\def\proof{{\indent\bf {\it Proof}}}
\def\endproof{\hfill\hbox{$\sqcup$}\llap{\hbox{$\sqcap$}}\medskip}

\makeatletter

\title[A lower bound for eigenvalues ]{ A lower bound for
eigenvalues of \\ the poly-Laplacian with arbitrary order* }
\author{Qing-Ming Cheng, Xuerong Qi\ and\ Guoxin Wei}
\address{Qing-Ming Cheng\\  Department of Mathematics, Faculty of Science and Engineering,
Saga University, Saga 840-8502,  Japan, cheng@ms.saga-u.ac.jp}
\address{Xuerong Qi\\  Department of Mathematics, Faculty of Science and Engineering,
Saga University, Saga 840-8502,  Japan, qixuerong609@gmail.com}
\address{Guoxin Wei\\  School of Mathematical Sciences, South China Normal University,
510631, Guangzhou,  China, weigx03@mails.tsinghua.edu.cn}

\begin{document}
\maketitle

\begin{abstract}
\noindent In this paper, we study eigenvalues of the poly-Laplacian with arbitrary order
on a bounded domain in an $n$-dimensional Euclidean space and obtain a lower bound for eigenvalues, 
which gives an important improvement of results due to Levine and Protter \cite{[LP]}. In particular, the result of Melas \cite{[M]} is included here.

\end{abstract}



\footnotetext{* Research partially supported by a Grant-in-Aid for
Scientific Research from JSPS and NSFC.}

\section *{1. Introduction}

Let $\Omega\subset\mathbb{R}^n$ be  a bounded domain with piecewise
smooth boundary $\partial\Omega$ in an $n$-dimensional Euclidean
space
 $\mathbb{R}^n$. Let $\lambda_i$ be the $i$-th eigenvalue of the Dirichlet eigenvalue
problem of the poly-Laplacian with arbitrary order:
\begin{equation*}
  {\begin{cases}
     (-\Delta)^l u =\lambda u,& \ \ {\rm in} \ \ \Omega ,\\
     u=\dfrac{\partial u}{\partial \nu}=\cdots=\dfrac{\partial^{l-1} u}{\partial \nu^{l-1}}=0, & \ \ {\rm on}  \ \ \partial \Omega,
     \end{cases}}
     \eqno{(1.1)}
\end{equation*}
where $\Delta$ is the Laplacian in $\mathbb{R}^n$ and $\nu$ denotes
the outward unit normal vector field of the boundary $\partial
\Omega$.
 It is well known that the spectrum of this eigenvalue problem
is real and discrete:
$$0<\lambda_1\leq\lambda_2\leq\lambda_3\leq\cdots\rightarrow+\infty,$$
where each $\lambda_i$ has finite multiplicity which is repeated
according to its multiplicity.

Let $V(\Omega)$ denote the volume of $\Omega$ and let $B_n$ denote
the volume of the unit ball in $\mathbb{R}^n$. When $l=1$, the
eigenvalue problem (1.1) is called a fixed membrane problem. In this
case, one has
 the following Weyl's asymptotic formula
$$\lambda_k\sim\frac{4\pi^2}{(B_nV(\Omega))^{\frac{2}{n}}}k^{\frac{2}{n}},
\ \ \ k\rightarrow+\infty.\eqno{(1.2)}$$ From the above asymptotic
formula, one can obtain
$$\frac{1}{k}\sum_{i=1}^k\lambda_i\sim\frac{n}{n+2}\frac{4\pi^2}{(B_nV(\Omega))^{\frac{2}{n}}}k^{\frac{2}{n}},
\ \ \ k\rightarrow+\infty.\eqno{(1.3)}$$
   P\'olya \cite{[P]} proved that
$$\lambda_k\geq\frac{4\pi^2}{(B_nV(\Omega))^{\frac{2}{n}}}k^{\frac{2}{n}},\
\ \ \rm{for}\  k=1,2,\cdots,\eqno{(1.4)}$$ if $\Omega$ is a tiling
domain in $\mathbb{R}^n$. Moreover,  he proposed  the following:

\par\bigskip \noindent{\bf Conjecture of P\'olya}. {\it
If $\Omega$ is a bounded domain in $\mathbb{R}^n$, then the $k$-th
eigenvalue $\lambda_k$ of the fixed membrane problem satisfies
$$\lambda_k\geq\frac{4\pi^2}{(B_nV(\Omega))^{\frac{2}{n}}}k^{\frac{2}{n}},\
\ \ {\rm for}\ k=1,2,\cdots.\eqno{(1.5)}$$} On the conjecture of
P\'olya, Berezin \cite{[B]} and Lieb  \cite{[L]} gave a partial
solution. In particular, Li and Yau \cite{[LY]} proved the following
$$
\frac{1}{k}\sum_{i=1}^k\lambda_i\geq\frac{n}{n+2}\frac{4\pi^2}{(B_nV(\Omega))^{\frac{2}{n}}}k^{\frac{2}{n}},\
\ \ \rm{for}\ k=1,2,\cdots.\eqno{(1.6)}
$$
The formula (1.3) shows that the result of Li and Yau is sharp in the sense of average.
 From this formula (1.6), one can derive
$$\lambda_k\geq\frac{n}{n+2}\frac{4\pi^2}{(B_nV(\Omega))^{\frac{2}{n}}}k^{\frac{2}{n}},\
\ \ \rm{for}\ k=1,2,\cdots,\eqno{(1.7)}$$ which gives a partial
solution for the conjecture of P\'olya with a factor
$\frac{n}{n+2}$.
 Recently, Melas \cite{[M]} have improved the estimate (1.6) to the following:
$$\frac{1}{k}\sum_{i=1}^k\lambda_i\geq\frac{n}{n+2}\frac{4\pi^2}{(B_nV(\Omega))^{\frac{2}{n}}}k^{\frac{2}{n}}
+\frac{1}{24(n+2)}\frac{V(\Omega)}{I(\Omega)},\ \ \ \rm{for}\
k=1,2,\cdots,\eqno{(1.8)}$$ where
$$
I(\Omega)=\min\limits_{a\in \mathbb{R}^n}\int_\Omega|x-a|^2dx
$$
 is called {\it  the moment of inertia} of
$\Omega$.

When $l=2$, the eigenvalue problem (1.1) is called a clamped plate problem.
 For the eigenvalues of the clamped plate problem, it follows from  Agmon  \cite{[Ag]} and Pleijel \cite{[P1]}
that
$$
\lambda_k\sim
\dfrac{16\pi^4}{\big(B_nV(\Omega)\big)^{\frac{4}{n}}}k^{\frac{4}{n}},\
\ \  k\rightarrow+\infty.
  \eqno{(1.9)}
  $$
This implies that
$$
\frac{1}{k}\sum_{i=1}^k\lambda_i
\sim\frac{n}{n+4}\dfrac{16\pi^4}{\big(B_nV(\Omega)\big)^{\frac{4}{n}}}k^{\frac{4}{n}},
\ \ k\rightarrow+\infty.\eqno{(1.10)}
$$
Furthermore, Levine and Protter \cite{[LP]} proved that the eigenvalues of
the clamped plate problem satisfy
$$\dfrac{1}{k}\sum_{i=1}^k\lambda_i
\geq\frac{n}{n+4}\dfrac{16\pi^4}{\big(B_nV(\Omega)\big)^{\frac{4}{n}}}k^{\frac{4}{n}}.\eqno{(1.11)}$$
The formula (1.10) shows that the coefficient
of $k^{\frac{4}{n}}$ is the best possible constant. Very recently, the first author and the third author \cite{[CW]}
obtained the following estimate which is an improvement of (1.11):
$$
\aligned \frac{1}{k}\sum_{i=1}^k\lambda_i&\geq
\frac{n}{n+4}\dfrac{16\pi^4}{\big(B_nV(\Omega)\big)^{\frac{4}{n}}}k^{\frac{4}{n}}\\
  & \ \ \ +c_n\frac{n}{n+2}
 \dfrac{4\pi^2}{\big(B_nV(\Omega)\big)^{\frac{2}{n}}}\frac{V(\Omega)}{I(\Omega)}k^{\frac{2}{n}}+d_n
 \left(\frac{V(\Omega)}{I(\Omega)}\right)^2,
\endaligned
\eqno{(1.12)}
$$
where  $c_n$ and $d_n$ are constants depending only on the dimension $n$.

When $l\geq 3$, Levine and Protter \cite{[LP]} proved the following
$$
\frac{1}{k}\sum_{i=1}^k\lambda_i\geq\frac{n}{n+2l}\frac{(2\pi)^{2l}}{(B_nV(\Omega))^{\frac{2l}{n}}}k^{\frac{2l}{n}},\
\ \ \rm{for}\ k=1,2,\cdots.\eqno{(1.13)}
$$
From the above formula, one can obtain
$$
\lambda_k\geq\frac{n}{n+2l}\frac{(2\pi)^{2l}}{(B_nV(\Omega))^{\frac{2l}{n}}}k^{\frac{2l}{n}},\
\ \ \rm{for}\ k=1,2,\cdots.\eqno{(1.14)}
$$

In this paper we investigate eigenvalues of the Dirichlet eigenvalue problem (1.1) of Laplacian with any order.
We give an important improvement of the result (1.13) due to Levine and Protter \cite{[LP]}
by adding  $l$ terms of lower order of $k^{\frac{2l}{n}}$ to its right hand side.  In fact, we prove the
following:

\par\bigskip \noindent{\bf Theorem 1}. {\it Let $\Omega$ be a
bounded domain in an $n$-dimensional Euclidean space $\mathbb{R}^n$.
Assume that $\lambda_i$, $i=1,2,\cdots,$ is the $i$-th eigenvalue of
the eigenvalue problem {\rm (1.1)}. Then the eigenvalues satisfy
$$
\aligned \frac{1}{k}\sum_{j=1}^k\lambda_j&\geq
\dfrac{n}{n+2l}\dfrac{(2\pi)^{2l}}{\big(B_nV(\Omega)\big)^{\frac{2l}{n}}}k^{\frac{2l}{n}}\\
 & +\frac{n}{n+2l}\sum_{p=1}^{l}
 \frac{(l+1-p)}{(24)^p n\cdots(n+2p-2)}\dfrac{(2\pi)^{2(l-p)}}{(B_nV(\Omega))^{\frac{2(l-p)}{n}}}
 \left(\frac{V(\Omega)}{I(\Omega)}\right)^pk^{\frac{2(l-p)}{n}}.
\endaligned
$$}

\par\bigskip \noindent{\bf Remark 1}. If we take $l=1$ in Theorem 1, then we obtain the inequality (1.8).

\section*{2. Proof of Theorem}

In this section, we firstly introduce some definitions and basic
facts about the symmetric decreasing rearrangements. Next, we give
the proof of Theorem 1.

For a bounded domain $\Omega\subset \mathbb{R}^n$, {\it the moment
of inertia} of $\Omega$ is defined by
$$
I(\Omega)=\min\limits_{a\in \mathbb{R}^n}\int_{\Omega}|x-a|^2dx.
$$
By translating the origin, we may assume that
$$
I(\Omega)=\int_\Omega|x|^2dx.
$$
Let $\Omega^*$ be the symmetric rearrangement of $\Omega$, that is, $\Omega^*$ is the open ball centered at the
origin with the same volume as
$\Omega$. Then
$$\Omega^*=\left\{x\in
\mathbb{R}^n; \  |x|<\biggl(\dfrac{{\rm V}
(\Omega)}{B_n}\biggl)^{\frac1n}\right\}.$$ By using the symmetric
rearrangement $\Omega^*$ of $\Omega$, we have
$$
I(\Omega)=\int_\Omega |x|^2dx\geq\int_{\Omega^*}
|x|^2dx=\frac{n}{n+2}V(\Omega)\left(\frac{V(\Omega)}{B_n}\right)^{\frac{2}{n}}.\eqno{(2.1)}
$$

 Let $f$ be a nonnegative continuous
function on $\Omega$. We consider its {\it distribution
function} $\mu_f(t)$ defined by
$$
\mu_f(t)={\rm Vol}(\{x\in\Omega; \ f(x)>t\}).
$$
The distribution function can be viewed as a function from
$[0,+\infty)$ to $[0, V(\Omega)]$. The {\it symmetric decreasing
rearrangement $f^*$ of $f$} is defined by
$$f^*(x)=\inf \{t\geq0; \ \mu_f(t)<B_n|x|^n\},\ \ \ {\rm for} \ x\in \Omega^*.$$ By definition, we know that
$f^*(x)$ is a radially symmetric function and $${\rm Vol}(\{x\in
\Omega; \ f(x)>t\})={\rm Vol}(\{x\in \Omega^*; \ f^*(x)>t\}), \ \ \forall \
t>0.$$

Let $f^*(x)=\phi(|x|)$. Then one gets that $\phi: [0, +\infty)\rightarrow
[0, \sup {f}]$ is a decreasing function of $|x|$. We may assume that
$\phi$ is absolutely continuous.
 It is well known that
$$\int_{\Omega}f(x)dx=\int_{\Omega^{*}}f^*(x)dx=nB_n\int_0^{+\infty}
s^{n-1}\phi(s)ds\eqno{(2.2)}
$$
and
$$\int_{\Omega}|x|^{2l}f(x)dx\geq\int_{\Omega^{*}}|x|^{2l}f^*(x)dx=nB_n\int_0^{+\infty}
s^{n+2l-1}\phi(s)ds.\eqno{(2.3)}
$$
Good sources of further information on rearrangements are
\cite{[Ba]}, \cite{[Po]}.

One gets from the coarea formula that
$$
\mu_f(t)=\int_t^{\sup f}\int_{\{f=s\}}|\nabla f|^{-1}d\sigma_sds.
$$
Since $f^*$ is radial, we have
\begin{equation*}
\begin{aligned}
\mu_f(\phi(s))&={\rm Vol}\{x\in
\Omega; \ f(x)>\phi(s)\}={\rm Vol}\{x\in \Omega^*; \ f^*(x)>\phi(s)\}\\
&={\rm Vol}\{x\in \Omega^*; \ \phi(|x|)>\phi(s)\}=B_ns^n.
\end{aligned}
\end{equation*}
 It follows that
$$
nB_ns^{n-1}=\mu_f^{'}(\phi(s))\phi^{'}(s)
$$
for almost every $s$. Putting $\tau:=\sup |\nabla f|$, we obtain
from the above equations and the isoperimetric inequality that
$$
\aligned -\mu_f^{'}(\phi(s))&=\int_{\{f=\phi(s)\}}|\nabla
   f|^{-1}d\sigma_{\phi(s)}\geq\tau^{-1}{\rm Vol}_{n-1}(\{f=\phi(s)\})\\
   &\geq\tau^{-1}nB_ns^{n-1}.
\endaligned
$$
Therefore, one  obtains
$$-\tau\leq \phi^{'}(s)\leq 0\eqno{(2.4)}$$
for almost every $s$.

In order to prove our theorem, we need the following lemma.

\par\bigskip \noindent{\bf Lemma 1}. {\it Let $b\geq1$, $\eta, A>0$
and $\psi: [0, +\infty)\rightarrow [0, +\infty)$ be a
decreasing, absolutely continuous function such that
$$-\eta\leq \psi^{'}(s)\leq 0,\ \ A=\int_0^{+\infty}s^{b-1}\psi(s)ds.$$
For any positive integer $l$, let
$$
A_{l}:=\int_0^{+\infty}s^{b+2l-1}\psi(s)ds.
$$
Then, we have
$$ A_{l}
 \geq\dfrac{1}{b+2l}\left[(bA)^{\frac{b+2l}{b}}\psi(0)^{-\frac{2l}{b}}
    +\sum_{p=1}^{l}\frac{(l+1-p)(bA)^{\frac{b+2(l-p)}{b}}
    \psi(0)^{\frac{2pb-2(l-p)}{b}}}{6^p b\cdots(b+2p-2)\eta^{2p}}\right].
$$}

\vskip 3pt \noindent  \proof. \ \  Using the method of induction.
Firstly, one can get from the lemma of \cite{[M]} that
$$
A_1=\int_0^{+\infty} s^{b+1}\psi(s)ds\geq\frac{1}{b+2}\left[(bA)^{\frac{b+2}{b}}\psi(0)^{-\frac{2}{b}}
+\frac{A\psi(0)^2}{6\eta^2}\right].\eqno{(2.5)}
$$
Secondly, we assume that
$$ A_{r}
 \geq\dfrac{1}{b+2r}\left[(bA)^{\frac{b+2r}{b}}\psi(0)^{-\frac{2r}{b}}
    +\sum_{p=1}^{r}\frac{(r+1-p)(bA)^{\frac{b+2(r-p)}{b}}
    \psi(0)^{\frac{2pb-2(r-p)}{b}}}{6^p b\cdots(b+2p-2)\eta^{2p}}\right].
$$
Since the formula (2.5) holds for any $b\geq 1$, we have
$$
\aligned &\ \ \ \ A_{r+1}=\int_0^{+\infty} s^{b+2r+1}\psi(s)ds\\
 &\geq\frac{1}{b+2r+2}\left\{[(b+2r)A_{r}]^{\frac{b+2r+2}{b+2r}}\psi(0)^{-\frac{2}{b+2r}}
 +\frac{A_{r}\psi(0)^2}{6\eta^2}\right\}\\
 &\geq\dfrac{\psi(0)^{-\frac{2}{b+2r}}}{b+2r+2}
 \left[(bA)^{\frac{b+2r}{b}}\psi(0)^{-\frac{2r}{b}}
    +\sum_{p=1}^{r}\frac{(r+1-p)(bA)^{\frac{b+2(r-p)}{b}}
    \psi(0)^{\frac{2pb-2(r-p)}{b}}}{6^p b\cdots(b+2p-2)\eta^{2p}}\right]^{\frac{b+2r+2}{b+2r}}\\
    & \ \ \ +\frac{1}{(b+2r)(b+2r+2)}\sum_{p=1}^{r}\frac{(r+1-p)(bA)^{\frac{b+2(r-p)}{b}}
    \psi(0)^{\frac{2(p+1)b-2(r-p)}{b}}}{6^{p+1} b\cdots(b+2p-2)\eta^{2p+2}}\\
    & \ \ \ +\frac{(bA)^{\frac{b+2r}{b}}\psi(0)^{\frac{2b-2r}{b}}}{6(b+2r)(b+2r+2)\eta^2}
    \endaligned$$
 $$\aligned
   &=\dfrac{\psi(0)^{-\frac{2}{b+2r}}}{b+2r+2}\left[(bA)^{\frac{b+2r}{b}}\psi(0)^{-\frac{2r}{b}}\right]^{\frac{b+2r+2}{b+2r}}\\
    & \ \ \ \times\left[1+\sum_{p=1}^{r}\frac{(r+1-p)(bA)^{\frac{-2p}{b}}
    \psi(0)^{\frac{2pb+2p}{b}}}{6^p b\cdots(b+2p-2)\eta^{2p}}\right]^{\frac{b+2r+2}{b+2r}}
    +\frac{(bA)^{\frac{b+2r}{b}}\psi(0)^{\frac{2b-2r}{b}}}{6(b+2r)(b+2r+2)\eta^2}\\
    & \ \ \ +\frac{1}{(b+2r)(b+2r+2)}\sum_{p=2}^{r+1}\frac{(r+2-p)(bA)^{\frac{b+2r-2p+2}{b}}
    \psi(0)^{\frac{2pb-2r+2p-2}{b}}}{6^{p} b\cdots(b+2p-4)\eta^{2p}}\\
   &=\dfrac{(bA)^{\frac{b+2r+2}{b}}\psi(0)^{-\frac{2r+2}{b}}}{b+2r+2}
   \left[1+\sum_{p=1}^{r}\frac{(r+1-p)(bA)^{\frac{-2p}{b}}
    \psi(0)^{\frac{2pb+2p}{b}}}{6^p b\cdots(b+2p-2)\eta^{2p}}\right]^{\frac{b+2r+2}{b+2r}}\\
    & \ \ \ +\frac{(bA)^{\frac{b+2r}{b}}\psi(0)^{\frac{2b-2r}{b}}}{6(b+2r)(b+2r+2)\eta^2}
       +\frac{(bA)
    \psi(0)^{2(r+1)}}{6^{r+1}b\cdots(b+2r+2)\eta^{2(r+1)}}\\
     & \ \ \ +\frac{1}{(b+2r)(b+2r+2)}\sum_{p=2}^{r}\frac{(r+2-p)(bA)^{\frac{b+2(r+1-p)}{b}}
    \psi(0)^{\frac{2pb-2(r+1-p)}{b}}}{6^{p} b\cdots(b+2p-4)\eta^{2p}}.
    \endaligned$$
It follows from the Taylor formula that
    $$\aligned
     A_{r+1}&\geq\dfrac{1}{b+2r+2}(bA)^{\frac{b+2r+2}{b}}\psi(0)^{-\frac{2r+2}{b}}\\
    & \ \ \ \times\left[1+\frac{b+2r+2}{b+2r}\sum_{p=1}^{r}\frac{(r+1-p)(bA)^{\frac{-2p}{b}}
    \psi(0)^{\frac{2pb+2p}{b}}}{6^p b\cdots(b+2p-2)\eta^{2p}}\right]\\
    & \ \ \ +\frac{(bA)^{\frac{b+2r}{b}}\psi(0)^{\frac{2b-2r}{b}}}{6(b+2r)(b+2r+2)\eta^2}
       +\frac{(bA)
    \psi(0)^{2(r+1)}}{6^{r+1}b\cdots(b+2r+2)\eta^{2(r+1)}}\\
     & \ \ \ +\frac{1}{(b+2r)(b+2r+2)}\sum_{p=2}^{r}\frac{(r+2-p)(bA)^{\frac{b+2(r+1-p)}{b}}
    \psi(0)^{\frac{2pb-2(r+1-p)}{b}}}{6^{p} b\cdots(b+2p-4)\eta^{2p}}\\
&=\frac{1}{b+2r+2}(bA)^{\frac{b+2r+2}{b}}\psi(0)^{-\frac{2r+2}{b}} \\
& \ \ \ +\frac{1}{b+2r}\sum_{p=1}^{r}\frac{(r+1-p)(bA)^{\frac{b+2(r+1-p)}{b}}
    \psi(0)^{\frac{2pb-2(r+1-p)}{b}}}{6^p b\cdots(b+2p-2)\eta^{2p}}\\
    & \ \ \ +\frac{(bA)^{\frac{b+2r}{b}}\psi(0)^{\frac{2b-2r}{b}}}{6(b+2r)(b+2r+2)\eta^2}
        +\frac{(bA)\psi(0)^{2(r+1)}}{6^{r+1} b\cdots(b+2r+2)\eta^{2(r+1)}}\\
    & \ \ \ +\frac{1}{(b+2r)(b+2r+2)}\sum_{p=2}^{r}\frac{(r+2-p)(bA)^{\frac{b+2(r+1-p)}{b}}
    \psi(0)^{\frac{2pb-2(r+1-p)}{b}}}{6^{p} b\cdots(b+2p-4)\eta^{2p}}\\
    &=\frac{1}{b+2r+2}(bA)^{\frac{b+2r+2}{b}}\psi(0)^{-\frac{2r+2}{b}}\\
    & \ \ \ +\left[\frac{r}{b(b+2r)}+\frac{1}{(b+2r)(b+2r+2)}\right]
    \frac{1}{6\eta^{2}}(bA)^{\frac{b+2r}{b}}\psi(0)^{\frac{2b-2r}{b}}
     \endaligned$$
$$\aligned
    & \ \ \ +\sum_{p=2}^{r}\left[\frac{r+1-p}{b+2r}+\frac{(r+2-p)(b+2p-2)}{(b+2r)(b+2r+2)}\right]
    \frac{(bA)^{\frac{b+2(r+1-p)}{b}}
    \psi(0)^{\frac{2pb-2(r+1-p)}{b}}}{6^p b\cdots(b+2p-2)\eta^{2p}}\\
     & \ \ \ +\frac{(bA)\psi(0)^{2(r+1)}}{6^{r+1} b\cdots(b+2r+2)\eta^{2(r+1)}}\\
    &\geq\dfrac{1}{b+2(r+1)}(bA)^{\frac{b+2(r+1)}{b}}\psi(0)^{-\frac{2(r+1)}{b}}\\
    & \ \ \ +\dfrac{1}{b+2(r+1)}\sum_{p=1}^{r+1}\frac{(r+2-p)(bA)^{\frac{b+2(r+1-p)}{b}}
    \psi(0)^{\frac{2pb-2(r+1-p)}{b}}}{6^p b\cdots(b+2p-2)\eta^{2p}}.
    \endaligned
$$
This completes the proof of Lemma 1.\endproof

\vskip 3pt\noindent  \proof \ {\it of  Theorem 1}.  Let $u_j$ be an orthonormal eigenfunction
corresponding to the eigenvalue $\lambda_j$, that is, $u_j$ satisfies
\begin{equation*}
  {\begin{cases}
     (-\Delta)^l u_j = \lambda_j  u_j,& \ \ {\rm in} \ \ \Omega ,\\
     u_j=\dfrac{\partial u_j}{\partial \nu}=\cdots=\dfrac{\partial^{l-1} u_{j}}{\partial \nu^{l-1}}=0 , & \ \ {\rm on}  \ \ \partial \Omega, \\
     \int_{\Omega} u_iu_j=\delta_{ij}, & \ \ \text{for any $i$, $j$}.
     \end{cases}}
     \eqno{(2.6)}
\end{equation*}
Thus, $\{u_j\}_{j=1}^{\infty}$ forms an orthonormal basis of $L^2(\Omega)$.
We define a function $\varphi_j$ by
\begin{equation*}
   \varphi_j(x)={\begin{cases}
    u_j(x), & \ \  x\in \Omega ,\\
     0 , & \ \ x\in \mathbb{R}^n\setminus\Omega. \\
         \end{cases}}
         \eqno{(2.7)}
\end{equation*}
The
Fourier transform $\widehat{\varphi}_j(z)$ of $\varphi_j(x)$ is then given by
$$\widehat{\varphi}_j(z)=(2\pi)^{-n/2}\int_{\mathbb{R}^n} \varphi_j(x)e^{i<x,z>}dx
=(2\pi)^{-n/2}\int_\Omega u_j(x)e^{i<x,z>}dx.\eqno{(2.8)}$$

We fix a $k\geq 1$ and set
$$f(z)=\sum_{j=1}^{k}|\widehat{\varphi}_j(z)|^{2}, \ \ \ \ {\rm for } \ z\in\mathbb{R}^n.$$
From
 Bessel's inequality, it follows that
$$\aligned 0\leq f(z)&=\sum_{j=1}^{k}|\widehat{\varphi}_j(z)|^{2}=(2\pi)^{-n}\sum_{j=1}^k
\left|\int_{\Omega}u_j(x)e^{i<x,z>}dx\right|^{2}\\
&\leq(2\pi)^{-n}\int_\Omega|e^{i<x,z>}|^2dx=(2\pi)^{-n}V(\Omega).\endaligned\eqno{(2.9)}$$
By Parseval's identity, we have
$$
\int_{\mathbb{R}^n}f(z)dz=\sum_{j=1}^k\int_{\mathbb{R}^n}
|\widehat{\varphi}_j(z)|^2dz=\sum_{j=1}^k\int_{\mathbb{R}^n}
\varphi_j^2(x)dx=\sum_{j=1}^k\int_{\Omega}
u_j^2(x)dx=k.\eqno{(2.10)}
$$
Furthermore, we deduce from integration by parts and Parseval's identity that
$$\aligned\int_{\mathbb{R}^n}|z|^{2l}f(z)dz&=\sum_{j=1}^k\int_{\mathbb{R}^n}|z|^{2l}|\widehat{\varphi}_j(z)|^2dz\\
&=\sum_{j=1}^k\int_{\mathbb{R}^n}|z|^{2l}\left|(2\pi)^{-n/2}\int_\Omega u_j(x)e^{i<x,z>}dx \right|^{2}dz\\
&=\sum_{j=1}^k\sum_{r_{1},\cdots,r_{l}=1}^{n}\int_{\mathbb{R}^n}
\left|(2\pi)^{-n/2}\int_\Omega z_{r_{1}}\cdots z_{r_{l}}u_j(x)e^{i<x,z>}dx \right|^{2}dz\\
&=\sum_{j=1}^k\sum_{r_{1},\cdots,r_{l}=1}^{n}\int_{\mathbb{R}^n}
\left|(2\pi)^{-n/2}\int_\Omega u_j(x)\frac{\partial^{l} e^{i<x,z>}}
{\partial x_{r_{1}}\cdots\partial x_{r_{l}}}dx \right|^{2}dz\\
&=\sum_{j=1}^k\sum_{r_{1},\cdots,r_{l}=1}^{n}\int_{\mathbb{R}^n}
\left|(2\pi)^{-n/2}\int_\Omega \frac{\partial^{l}u_j(x)}{\partial x_{r_{1}}\cdots\partial x_{r_{l}}}e^{i<x,z>}dx \right|^{2}dz\\
&=\sum_{j=1}^k\sum_{r_{1},\cdots,r_{l}=1}^{n}\int_{\mathbb{R}^n}\left|\widehat{\frac{\partial^{l}u_j}{\partial x_{r_{1}}\cdots\partial x_{r_{l}}}}\right|^{2}dz\\
&=\sum_{j=1}^k\sum_{r_{1},\cdots,r_{l}=1}^{n}\int_{\mathbb{R}^n}\left(\frac{\partial^{l}u_j}{\partial x_{r_{1}}\cdots\partial x_{r_{l}}}\right)^{2}dx\\
&=\sum_{j=1}^k\int_\Omega u_j(-\Delta)^lu_jdx=\sum_{j=1}^k\lambda_j.\\
\endaligned\eqno{(2.11)}$$
Since
$$\nabla \widehat{\varphi}_j(z)=(2\pi)^{-n/2}\int_\Omega ixu_j(x)e^{i<x,z>}dx,\eqno{(2.12)}$$
  we obtain from Bessel's inequality
$$\sum_{j=1}^k|\nabla \widehat{\varphi}_j(z)|^2\leq
(2\pi)^{-n}\int_\Omega|ixe^{i<x,z>}|^2dx=(2\pi)^{-n}I(\Omega).\eqno{(2.13)}$$
It follows from (2.9), (2.13) and the Cauchy-Schwarz inequality that
$$
\aligned |\nabla f(z)|&\leq
2\left(\sum_{j=1}^k|\widehat{\varphi}_j(z)|^2\right)^{1/2}\left(\sum_{j=1}^k|\nabla
     \widehat{\varphi}_j(z)|^2\right)^{1/2}\\
     &\leq2(2\pi)^{-n}\sqrt{V(\Omega)I(\Omega)}
\endaligned
\eqno{(2.14)}
$$
for every $z\in \mathbb{R}^n$.

Using the symmetric decreasing rearrangement $f^{*}$ of $f$ and noting
that
$$f^{*}(x)=\phi(|x|),\ \ \tau=\sup |\nabla
f|\leq2(2\pi)^{-n}\sqrt{V(\Omega)I(\Omega)}:=\eta,$$ we obtain, from
(2.4),
$$
-\eta\leq-\tau\leq \phi^{'}(s)\leq 0\eqno{(2.16)}
$$
for almost every $s$. According to (2.2) and (2.10), we infer
$$k=\int_{\mathbb{R}^n}f(z)dz=\int_{\mathbb{R}^n}f^*(z)dz=nB_n\int_0^{+\infty}
s^{n-1}\phi(s)ds.\eqno{(2.17)}
$$
From (2.3) and (2.11), we obtain
$$
\sum_{j=1}^k\lambda_j=\int_{\mathbb{R}^n}|z|^{2l}f(z)dz\geq\int_{\mathbb{R}^n}|z|^{2l}f^*(z)dz
     =nB_n\int_0^{+\infty} s^{n+2l-1}\phi(s)ds.
\eqno{(2.18)}
$$

Now, we can apply Lemma 1 to the function $\phi$ with
$$ b=n, \ \ \ \ A=\frac{k}{nB_n}, \ \ \ \ \eta=2(2\pi)^{-n}\sqrt{V(\Omega)I(\Omega)}.\eqno{(2.19)}$$
 We conclude that
 $$\aligned \sum_{j=1}^k\lambda_j&
\geq\dfrac{nB_n}{n+2l}\left(\frac{k}{B_n}\right)^{\frac{n+2l}{n}}\phi(0)^{-\frac{2l}{n}}\\
    & +\dfrac{nB_n}{n+2l}\sum_{p=1}^l\frac{(l+1-p)}{6^{p}n\cdots(n+2p-2)\eta^{2p}}\left(\frac{k}{B_n}\right)^{\frac{n+2l-2p}{n}}
    \phi(0)^{\frac{2pn+2p-2l}{n}}.
\endaligned
\eqno{(2.20)}
$$

Note that
 $0<\phi(0)\leq\sup f\leq(2\pi)^{-n}V(\Omega)$. Hence we consider the function
  $F$ defined by
$$\aligned F(t)&
=\dfrac{nB_n}{n+2l}\left(\frac{k}{B_n}\right)^{\frac{n+2l}{n}}t^{-\frac{2l}{n}}\\
& +\dfrac{nB_n}{n+2l}
\sum_{p=1}^l\frac{(l+1-p)}{6^{p}n\cdots(n+2p-2)\eta^{2p}}\left(\frac{k}{B_n}\right)^{\frac{n+2l-2p}{n}}
    t^{\frac{2pn+2p-2l}{n}},\endaligned
\eqno{(2.21)}$$ for $t\in(0, (2\pi)^{-n}V(\Omega)]$. From (2.1), we
have
$$\eta\geq(2\pi)^{-n}B_n^{-\frac{1}{n}}V(\Omega)^{\frac{n+1}{n}}.\eqno{(2.22)}$$
By a direct calculation, one gets from
$B_n=\dfrac{2\pi^{\frac{n}{2}}}{n\Gamma(\frac{n}{2})}$ that
$$\frac{B_n^{\frac{4}{n}}}{(2\pi)^{2}}<\frac{1}{2},\eqno{(2.23)}$$
where $\Gamma(\frac{n}{2})$ is the Gamma function. Thus, it follows
from (2.22) and (2.23) that
$$\aligned
 &  F^{'}(t)
=\dfrac{2B_nt^{-\frac{n+2l}{n}}}{n+2l}\left(\frac{k}{B_n}\right)^{\frac{n+2l}{n}}
\left[-l+\sum_{p=1}^l\frac{(l+1-p)(pn+p-l)
t^{\frac{2p(n+1)}{n}}}{6^{p}n\cdots(n+2p-2)\eta^{2p}}
\left(\frac{k}{B_n}\right)^{-\frac{2p}{n}}\right]\\
&\leq\dfrac{2B_n}{n+2l}\left(\frac{k}{B_n}\right)^{\frac{n+2l}{n}}t^{-\frac{n+2l}{n}}
\left[-l+\sum_{p>\frac{l}{n+1}}^l\frac{(l+1-p)(pn+p-l)}{6^{p}n\cdots(n+2p-2)}\left(\frac{B_n^{\frac{4}{n}}}{(2\pi)^{2}}\right)^{p}\right]\\
&<\dfrac{2B_n}{n+2l}\left(\frac{k}{B_n}\right)^{\frac{n+2l}{n}}t^{-\frac{n+2l}{n}}
\left[-l+\sum_{p>\frac{l}{n+1}}^l\frac{(l+1-p)(pn+p-l)}{(12)^{p}n\cdots(n+2p-2)}\right]\\
&<\dfrac{2B_n}{n+2l}\left(\frac{k}{B_n}\right)^{\frac{n+2l}{n}}t^{-\frac{n+2l}{n}}
\left[-l+\frac{l(n+1-l)}{12n}+\sum_{p>\frac{l}{n+1},p\neq
1}^l\frac{p^{2}n(n+1)}{(12)^{p}n\cdots(n+2p-2)}\right]\\
&<\dfrac{2B_n}{n+2l}\left(\frac{k}{B_n}\right)^{\frac{n+2l}{n}}t^{-\frac{n+2l}{n}}
\left[-l+\frac{l}{12}+\sum_{p>\frac{l}{n+1},p\neq
1}^l\frac{p^{2}}{(12)^{p}}\right]\\
&<\dfrac{2B_n}{n+2l}\left(\frac{k}{B_n}\right)^{\frac{n+2l}{n}}t^{-\frac{n+2l}{n}}
\left[-l+\frac{l}{12}+\frac{1}{12}\right]<0.
\endaligned
$$
We obtain that $F(t)$ is a decreasing function on $(0,
(2\pi)^{-n}V(\Omega)]$. Then we can replace $\phi(0)$ by
$(2\pi)^{-n}V(\Omega)$ in (2.20), namely
 $$\aligned \sum_{j=1}^k\lambda_j&
\geq\dfrac{n}{n+2l}\frac{(2\pi)^{2l}}{(B_nV(\Omega))^{\frac{2l}{n}}}k^{\frac{n+2l}{n}}\\
    & \ \ \ +\dfrac{n}{n+2l}\sum_{p=1}^l\frac{(l+1-p)}{6^{p}n\cdots(n+2p-2)\eta^{2p}}\frac{(V(\Omega))^{\frac{2pn+2p-2l}{n}}}
    {(2\pi)^{2pn+2p-2n}B_n^{\frac{2l-2p}{n}}}k^{\frac{n+2l-2p}{n}}\\
    &=\dfrac{n}{n+2l}\frac{(2\pi)^{2l}}{(B_nV(\Omega))^{\frac{2l}{n}}}k^{\frac{n+2l}{n}}\\
    & \ \ \ +\dfrac{n}{n+2l}\sum_{p=1}^l\frac{(l+1-p)}{24^{p}n\cdots(n+2p-2)}\frac{(2\pi)^{2(l-p)}}
    {(B_nV(\Omega))^{\frac{2(l-p)}{n}}}\left(\frac{V(\Omega)}{I(\Omega)}\right)^{p}k^{\frac{n+2(l-p)}{n}}.
\endaligned
$$
This completes the proof of Theorem 1.
\endproof

\end {document}